%% file: EJC_main.tex
\newif\ifarxiv

\arxivfalse

\ifarxiv
\documentclass[journal]{IEEEtran}
\else
 \documentclass[letterpaper, 10 pt, conference]{ieeeconf}
\fi

\usepackage{centernot}
\usepackage{cite}
\usepackage{amsmath,amssymb,amsfonts}

\usepackage{amsthm}

\usepackage{mathdots}
\usepackage{textcomp}
\usepackage[colorlinks=True,allcolors=blue]{hyperref}
\usepackage[margin=1in]{geometry}
\input{preamble.tex}
\usepackage[many]{tcolorbox}
\tcbuselibrary{breakable} 

\usepackage{soul}

\newcommand\QP{\hyperref[QP]{$\textbf{QP}$}\xspace}
\newcommand\QPslack{\hyperref[QPslack]{$\textbf{QPslack}$}\xspace}
\newcommand\QPslacktwo{\hyperref[QPslacktwo]{$\textbf{QPslack2}$}\xspace}
\newcommand\QPtilde{\hyperref[QPtilde]{$\widetilde{\textbf{QP}}$}\xspace}

\begin{document}
\title{
Solving Quadratic Programs with Slack Variables via ADMM \\ without Increasing the Problem Size} 

\author{Thomas Lew, Marcus Greiff, John Subosits, Brian Plancher
\thanks{T. Lew, M. Greiff, and J. Subosits are with Toyota Research Institute, Los Altos, USA. %
B. Plancher is with Barnard College, Columbia University and Dartmouth College, USA.}
\vspace{-2mm}
}
\maketitle

\begin{abstract}
Proximal methods such as the Alternating Direction Method of Multipliers (ADMM) are effective at solving constrained quadratic programs (QPs).  
To tackle infeasible QPs, slack variables are often introduced to ensure feasibility, which changes the structure of the problem, increases its size, and slows down numerical resolution. 
In this letter, we propose a simple ADMM scheme to tackle QPs with slack variables without increasing the size of the original problem. 
The only modification is a slightly different projection in the $z$-update, while the rest of the algorithm remains standard. 
We prove that the method is equivalent to applying ADMM to the QP with additional slack variables, even though slack variables are not added. 
Numerical experiments show speedups of the approach.
\end{abstract}

\section{Introduction}

Numerical optimization underpins a wide range of algorithms and applications in control. For example, model predictive control (MPC) algorithms steer systems by repeatedly solving optimization problems to select control inputs~\cite{Houska2018}, with many applications such as in  robotics~\cite{tranzatto2022cerberus,wensing2023optimization,adabag2024mpcgpu}, automotive~\cite{del2010automotive,LewGreiffICRA2025}, and power systems~\cite{venkat2008distributed,vazquez2014model}. 
MPC often relies on repeated solutions of constrained quadratic programs (QP) of the form
\begin{align}\label{QP} 
\tag{\QP}
\min_x\    
\tfrac{1}{2}x^\top Q x + q^\top x 
\ \ \text{s.t.}\ \, 
\ell \leq A x\leq u,
\end{align}
where $x\in\R^n$ are optimization variables, $Q\in\R^{n\times n}$ is positive semidefinite, $q\hspace{-1pt}\in\hspace{-1pt}\R^n$, $A\hspace{-1pt}\in\hspace{-1pt}\R^{p\times n}$, and $\ell\hspace{-1pt}\leq u\hspace{-1pt}\in\hspace{-1pt}\R^p$. Proximal methods such as the Alternating Direction Method of Multipliers (ADMM)~\cite{Boyd2011,Stellato2020,o2016conic} have proven effective in
solving such problems. %

In many applications, \QP might become infeasible. For example, a race car  
operating at the limits of friction~\cite{LewGreiffICRA2025} using MPC may transiently drive over road edges due to model mismatch. To maintain system operation in such situations, a solution that approximately solves \QP may be desired. 
For this reason, practitioners often solve QPs with soft constraints, where $\xi\in\R^p$ are additional slack variables and $\alpha>0$ is a penalization weight:
\begin{align}\label{QPslack} 
\tag{\QPslack}
\min_{x,\xi}\ \   
&\tfrac{1}{2}x^\top Q x + q^\top x  + \tfrac{\alpha}{2}\|\xi\|_2^2
\\
\text{s.t.}\ \ 
&\ell \leq A x+\xi\leq u.
\nonumber
\end{align}
By introducing the slack variables $\xi$ in the constraints, \QPslack is always feasible, and its solutions can be made as close as desired to those to \QP by selecting a large enough penalization weight $\alpha$. 
Considering the relaxed problem {\QPslack has a similar feasibility-enabling  effect as  penalizing constraints in the cost \cite{Feller2017,Oancea2023}.

\begin{figure}[!t]
\centering
\includegraphics[width=0.6\linewidth]{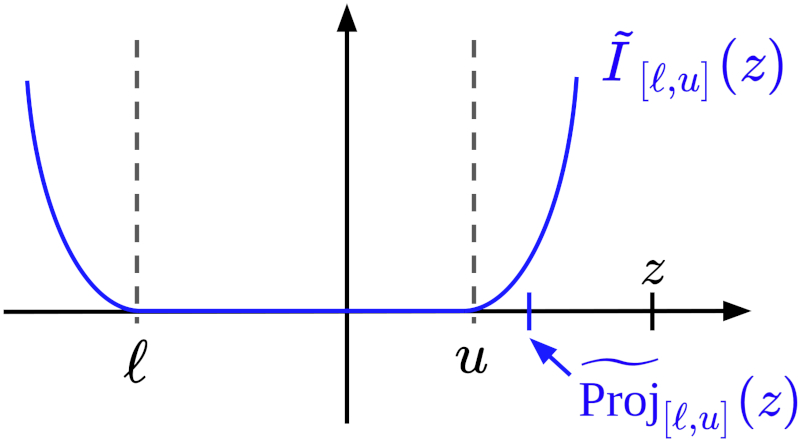}
\caption{Smoothed indicator function and projection operator.}
\vspace{-2mm}
    \label{fig:projection}
\end{figure} 

\textbf{Standard approach}: \QPslack   can  be reformulated as \QP, since it is equivalent to the following problem:
\begin{align}\label{QPslacktwo} 
\tag{\QPslacktwo}
\min_{\bar{x}=(x,\xi)}\  
&\tfrac{1}{2}\bar{x}^\top \bar{Q}
\bar{x} + \bar{q}^\top \bar{x}
\ \ \text{s.t.}\ \, 
\ell \leq \bar{A}\bar{x}\leq u,
\end{align}
where
$$
\bar{x}=\begin{bmatrix}x\\\xi
\end{bmatrix},
\ \ 
\bar{Q}=\begin{bmatrix}
Q&0\\0&\alpha
\end{bmatrix},
\ \ 
\bar{q}=\begin{bmatrix}q\\0
\end{bmatrix},
\ \ 
\bar{A}=\begin{bmatrix}A&I
\end{bmatrix}.
$$
Thus, \QPslack takes the form of \QP with additional optimization variables $\xi$, and it can be solved using a general-purpose \QP solver such as \texttt{OSQP} \cite{Stellato2020}. 
However, this approach has two limitations. First, it leads to slower numerical resolution compared to solving \QP (assuming that \QP is feasible), due to the %
increased problem size. Second, application-specific algorithms (e.g., for MPC) leverage the structure of $(Q,q,A,\ell,u)$ for efficiency (e.g.,~\cite{adabag2024mpcgpu}), which may change when introducing the slack variables $\xi$ via this standard reformulation, slowing down algorithm development and implementation.

\textbf{Contribution}: 
We present an ADMM approach for solving \QPslack that is faster and easier to implement than the standard approach. 
It consists of a slight modification to the standard ADMM projection step for \QP (see Table \ref{tab:admm}) that enables solving  \QPslack \emph{without} introducing slack variables, maintaining problem size and resulting in faster solves, while keeping the convergence guarantees of other ADMM methods. We give two derivations of the algorithm, one of which shows that the proposed method is just an ADMM scheme for \QPslack derived with an unusual splitting scheme that enables eliminating the slack variables.  
We release code of the ADMM approach and experiments open-source.  

\input{table}

\section{ADMM Background}
\subsection{ADMM for \QP}\label{sec:admm:qp}
We begin by describing a standard ADMM scheme for solving \QP as described in \cite{Boyd2011}. 
Throughout the manuscript, we rely on the following assumption on the matrices $(Q,A)$ defining \QP:

\begin{assumption}\label{assumption}
$Q+\rho A^\top A$ is invertible for any $\rho> 0$.
\end{assumption}

While we leave relaxing this assumption, e.g., using the splitting scheme in \cite{Stellato2020}, to future work, we note that QPs that do not satisfy Assumption \ref{assumption} can be slightly modified by adding sufficiently many inactive constraints $-\infty\leq x_i\leq \infty$ until $A^\top A$ is invertible.

We introduce auxiliary optimization variables $z\in\R^p$ to define the following equivalent problem to \QP:
\begin{align}\label{eq:QP:z}
\min_{x,z}\     
\tfrac{1}{2}x^\top Q x + q^\top x 
\ \ \text{s.t.}\ \, 
Ax=z, \ \  
\ell \leq z\leq u.
\end{align}
This problem is also equivalent to
\begin{align}
\label{eq:QP:indicator}
\min_{x,z}\    
&\tfrac{1}{2}x^\top Q x + q^\top x 
+ I_{[\ell, u]}(z) \ \  
\text{s.t.}\ \ 
Ax=z,
\end{align}
where $I_{[\ell,u]}$ is the indicator function
\begin{align*}
\nonumber
\\[-6mm]
I_{[\ell,u]}(z)\,{=}\begin{cases}
0&\hspace{-2mm}\text{if }\ell\leq z \leq u,
\\[-1mm]
\infty&\hspace{-2mm}\text{otherwise.}
\end{cases}
\end{align*}
We then introduce a step size parameter $\rho> 0$, and define a third equivalent optimization problem to \QP:
\begin{subequations}
\label{eq:QP:rho}
\begin{align}
\min_{x,z}\   
\     
&\tfrac{1}{2}x^\top Q x + q^\top x 
+ I_{[\ell, u]}(z) + \tfrac{\rho}{2}\|Ax-z\|^2
\\
\text{s.t.}\ \ \ 
&Ax=z.
\end{align}
\end{subequations}
Using the KKT multipliers $\mu\in\R^p$ associated to the constraints $Ax=z$, we define the augmented Lagrangian
\begin{align}\nonumber
\mathcal{L}(x,z,\mu)=\,
&\tfrac{1}{2}x^\top Q x + q^\top x 
+ I_{[\ell, u]}(z)+\tfrac{\rho}{2}\|Ax-z\|^2  %
\\
&\qquad+\mu^\top(Ax-z).\label{eq:QP:admm}
\end{align}
ADMM solves \QP by iteratively solving~\eqref{eq:QP:admm} until convergence through the following three-step process:
\begin{enumerate}[leftmargin=5mm]
\item \textbf{$x$-update}: minimize  \eqref{eq:QP:admm} over $x$ with $(z,\mu)$ fixed:
\begin{align}\label{eq:QP:admm:x}
x=(Q+\rho A^\top A)^{-1}(-q-A^\top\mu+\rho A^\top z).
\end{align}
\item \textbf{$z$-update}: minimize  \eqref{eq:QP:admm} over $z$ with $(x,\mu)$ fixed:
\begin{align}\label{eq:QP:admm:z}
&z=\text{Proj}_{[\ell,u]}(Ax+\mu/\rho),
\end{align}
where $\text{Proj}_{[\ell,u]}(\cdot)$ is the projection operator
\begin{align}\label{eq:proj}
\text{Proj}_{[\ell,u]}(\tilde{z})
=
\begin{cases}
\tilde{z}& \text{if }\ell \leq \tilde{z} \leq u,
\\[-1mm]
\ell
& \text{if }\tilde{z}< \ell,
\\[-1mm]
u
& \text{if }\tilde{z}> u.
\end{cases}
\end{align}
\item \textbf{dual update}: update the KKT multipliers via:
\begin{align}\label{eq:QP:admm:mu}
&\mu=\mu + \rho(Ax-z).
\end{align}
\end{enumerate}
The convergence of this standard ADMM scheme is well-studied, see  \cite{Boyd2011} and the references therein.

\textbf{Exit condition}. The convergence of the solver is assessed by the magnitude of the primal and dual residuals
\begin{equation}\label{eq:residuals}
r_{\textrm{prim}} = \|Ax-z\|,
\quad
r_{\textrm{dual}} = \|Qx+q+A^\top\mu\|.
\end{equation}
The solver returns a solution once $\max(r_{\textrm{prim}},r_{\textrm{dual}})< \epsilon$.

\textbf{Selecting $\rho$}. The step-size parameter $\rho$ plays an important role in the convergence speed of the solver \cite{Boyd2011,Stellato2020}. Every $n_{\rho}=25$ ADMM iterations,  we update $\rho$ as
\begin{equation}\label{eq:rho}
\rho\gets \begin{cases}
\rho\sqrt{\frac{r_{\textrm{prim}}}{ r_{\textrm{dual}}}} 
& \text{if $\sqrt{\frac{r_{\textrm{prim}}}{ r_{\textrm{dual}}}} > \kappa$ or $\sqrt{\frac{r_{\textrm{dual}}}{ r_{\textrm{prim}}}} > \kappa$,}
\\
\rho&\text{otherwise.}
\end{cases}
\end{equation}
As in \cite{Stellato2020}, we set the ratio parameter to $\kappa=5$.

\subsection{Standard ADMM for \QPslack}\label{sec:admm:qpslack}
Since \QPslack $=$ \QPslacktwo, and \QPslacktwo is an instance of \QP, we can use the standard ADMM scheme in the previous section to solve \QPslack.  
This ADMM scheme considers the following augmented Lagrangian 
\begin{align}\nonumber
\bar{\mathcal{L}}(\bar{x},z,\mu):=\,
&\tfrac{1}{2}\bar{x}^\top \bar{Q} \bar{x} + \bar{q}^\top \bar{x} 
+ I_{[\ell, u]}(z)+\frac{\rho}{2}\|\bar{A}\bar{x}-z\|^2  
\\[-1mm]
&\qquad+\mu^\top(\bar{A}\bar{x}-z),\label{eq:QP:admm:slack}
\end{align}
with $\bar{x}=(x,\xi)$, and consists of the following three steps:
\begin{enumerate}
\item \textbf{$(x,\xi)$-update}: minimize  \eqref{eq:QP:admm:slack} over $\bar{x}=(x,\xi)$ again with $(z,\mu)$ fixed:
\begin{align}\label{eq:QP:admm:x:slack}
\hspace{-4mm}&\bar{x}=(\bar{Q}+\rho \bar{A}^\top \bar{A})^{-1}(-\bar{q}-\bar{A}^\top \mu+\rho \bar{A}^\top z)
\\
\hspace{-4mm}&=
\begin{bmatrix}Q+\rho A^\top A & \rho A^\top
\\
\rho A & (\rho+\alpha)I\end{bmatrix}^{-1}
\begin{bmatrix}
-q-A^\top \mu+\rho A^\top z
\\
-\mu+\rho z
\end{bmatrix}.
\nonumber
\end{align}
\item \textbf{$z$-update}: minimize \eqref{eq:QP:admm:slack} over $z$ with $(\bar{x},\mu)$ fixed:
\begin{align}\label{eq:QP:admm:z:slack}
&z=\text{Proj}_{[\ell,u]}(\bar{A}\bar{x}+\mu/\rho).
\end{align}
\item \textbf{dual update}: update the KKT multipliers via
\begin{align}\label{eq:QP:admm:mu:slack}
&\mu=\mu + \rho(\bar{A}\bar{x}-z).
\end{align}
\end{enumerate}
By comparing the $x$-updates \eqref{eq:QP:admm:x} and \eqref{eq:QP:admm:x:slack}, we see that the $(x,\xi)$-update in \eqref{eq:QP:admm:x:slack} involves solving a larger linear system than for the problem \QP without the slack variables $\xi$, which can result in slower numerical resolution. %

The primal update in  \eqref{eq:QP:admm:x} is well-posed if the matrix $(\bar{Q}+\rho\bar{A}^\top\bar{A})$ is invertible, which holds true for sufficiently large $\alpha$ as a function of $(Q,A,\rho)$. %
Our proposed ADMM approach does not rely on this condition.

\section{ADMMSlack}\label{sec:admm:new}
We propose a simple ADMM scheme (\texttt{ADMMSlack}) to solve the QP with slack variables \QPslack \emph{without} adding the slack variables $\xi$ to the problem. In this section, we first derive the scheme via smoothed indicator functions. We will show in the next section that this scheme actually solves \QPslack.  
Consider the problem
\begin{align}\label{QPtilde} 
\tag{\QPtilde}
&\min_{x,z}\    
\tfrac{1}{2}x^\top Q x + q^\top x + \widetilde{I}_{[\ell, u]}(z)
\ \ 
\text{s.t.}\ \, 
A x=z,
\end{align}
where $\widetilde{I}_{[\ell, u]}(\cdot)$ denotes the smoothed indicator function defined as 
 $\widetilde{I}_{[\ell,u]}(z)=\sum_{i=1}^p
\widetilde{I}_i(z_i)$ with
$$
\widetilde{I}_i(z_i)=
\begin{cases}
0 & \text{if }\ell_i \leq z_i \leq u_i,
\\
\frac{\alpha}{2}(z_i-\ell_i)^2 
& \text{if }z_i<\ell_i,
\\
\frac{\alpha}{2}(z_i-u_i)^2 
& \text{if }z_i>u_i.
\end{cases}
$$
\QPtilde is similar to the reformulation of \QP in \eqref{eq:QP:indicator}, except that it leverages the smoothed indicator function $\widetilde{I}_{[\ell, u]}(\cdot)$ shown in Figure~\ref{fig:projection}. 
Note that \QPtilde is not a QP, and  that $\widetilde{I}_{[\ell, u]}$ approaches $I_{[\ell, u]}$ as $\alpha$ increases, so that $\widetilde{I}_{[\ell, u]}$  smoothly approximates   $I_{[\ell, u]}$.

Using the KKT multipliers $\mu\in\R^p$ associated with the constraints $Ax=z$, the proposed ADMM scheme considers the augmented Lagrangian
\begin{align}\nonumber
\widetilde{\mathcal{L}}(x,z,\mu):=\,
&\tfrac{1}{2}x^\top Q x + q^\top x 
+ \widetilde{I}_{[\ell, u]}(z)+\tfrac{\rho}{2}\|Ax-z\|^2  
\\
&\qquad+\mu^\top(Ax-z).\label{eq:QP:admm:new}
\end{align}

If we apply ADMM to \QPtilde and  sequentially minimize \eqref{eq:QP:admm:new}, we end up with the following three step processs: %
\begin{enumerate}[leftmargin=5mm]
\item \textbf{$x$-update}: minimize~\eqref{eq:QP:admm:new} over $x$ with $(z,\mu)$ fixed:
\begin{align*}
x=(Q+\rho A^\top A)^{-1}(-q-A^\top\mu+\rho A^\top z)=\eqref{eq:QP:admm:x}.
\end{align*}
\item \textbf{$z$-update}: minimize~\eqref{eq:QP:admm:new} over $z$ with $(x,\mu)$ fixed:
\begin{align}\label{eq:QP:admm:z:new}
&z =\widetilde{\text{Proj}}_{[\ell,u]}(Ax+\mu/\rho),
\end{align}
where $\widetilde{\text{Proj}}$ is the smoothed projection operator defined elementwise as
\begin{align}\label{eq:proj:tilde}
\widetilde{\text{Proj}}_{[\ell,u]}(\tilde{z})
=
\begin{cases}
\tilde{z}& \text{if }\ell \leq \tilde{z} \leq u,
\\
\frac{\rho \tilde{z} + \alpha \ell}{\rho+\alpha}
& \text{if }\tilde{z}< \ell,
\\
\frac{\rho \tilde{z} + \alpha u}{\rho+\alpha} 
& \text{if }\tilde{z}> u.
\end{cases}
\end{align}
\item \textbf{dual update}: update the KKT multipliers via
\begin{align*}
&\mu=\mu + \rho(Ax-z)=\eqref{eq:QP:admm:mu}.
\end{align*}
\end{enumerate}
This proposed ADMM scheme is more computationally and memory efficient than the standard ADMM scheme  for solving \QPslack (Section \ref{sec:admm:qpslack}), thanks to the $x$-update. The only difference with the ADMM scheme for solving \QP without slack variables $\xi$ is using a smoothed projection operator $\widetilde{\text{Proj}}_{[\ell,u]}$. Thus, this proposed scheme is easy to implement given an implementation of ADMM for solving \QP without slack variables $\xi$. 

\section{ADMMSlack is just an ADMM for \QPslack with a different splitting scheme}\label{sec:admm:proof}
\texttt{ADMMSlack} was derived as an ADMM on \QPtilde, with no clear connection to \QPslack. In this section, we show that \texttt{ADMMSlack} can be equivalently derived as an ADMM for solving \QPslack, where one sequentially optimizes $x$ and then $(z,\xi)$,  treating the slack variable $\xi$ as another auxiliary variable. 

First, note that \QPslack is also equivalent to
\begin{align*}  
\min_{x,z,\xi}\ \   
&\tfrac{1}{2}x^\top Q x + q^\top x  + \tfrac{\alpha}{2}\|\xi\|_2^2
\\
\text{s.t.}\ \ \, 
&A x=z, \  
\ell \leq z+\xi\leq u,
\end{align*}
which is equivalent to
\begin{align*}
\min_{x,z,\xi}\    
&\tfrac{1}{2}x^\top Q x + q^\top x   + \tfrac{\alpha}{2}\|\xi\|_2^2+ I_{[\ell, u]}(z+\xi)
\\
\text{s.t.}\ \, 
&A x=z.
\end{align*}
Then, we define the augmented Lagrangian
\begin{align}\nonumber
\mathcal{L}(x,z,\xi,\mu):=\,
&\tfrac{1}{2}x^\top Q x + q^\top x 
+ \tfrac{\alpha}{2}\|\xi\|_2^2+I_{[\ell, u]}(z+\xi) 
\\
&\ \ +\tfrac{\rho}{2}\|Ax-z\|^2+\mu^\top(Ax-z).\label{eq:QP:admm:slack:2}
\end{align}
Solving~\eqref{eq:QP:admm:slack:2} with ADMM then has the following three iterative steps:
\begin{enumerate}
\item \textbf{$x$-update}: minimize \eqref{eq:QP:admm:slack:2} over $x$ with $(z,\xi,\mu)$ fixed:
\begin{align*}
x&=(Q+\rho A^\top A)^{-1}(-q-A^\top\mu+\rho A^\top z)=\eqref{eq:QP:admm:x}.
\end{align*}
\item \textbf{$(z,\xi)$-update}: minimize  \eqref{eq:QP:admm:slack:2} over $(z,\xi)$ with $(x,\mu)$ fixed:
\begin{align}
\nonumber 
\min_{z,\xi}
\ \ 
&\frac{1}{2}\begin{bmatrix}z\\\xi\end{bmatrix}^\top
\begin{bmatrix}
\rho&0\\0&\alpha
\end{bmatrix}\begin{bmatrix}z\\\xi\end{bmatrix}
+
\begin{bmatrix}
-\rho Ax-\mu\\0
\end{bmatrix}^\top
\begin{bmatrix}z\\\xi\end{bmatrix}
\\
\text{s.t.}\ \ & \ell\leq z+\xi\leq u.
\label{eq:QP:zxi}
\end{align}
\item \textbf{dual update}: update the KKT multipliers via
\begin{align*}
&\mu=\mu + \rho(Ax-z)=\eqref{eq:QP:admm:mu}.
\end{align*}
\end{enumerate}
Note that the $x$ and dual updates are exactly \eqref{eq:QP:admm:x} and \eqref{eq:QP:admm:mu} of the original ADMM scheme for \QP. 
Moreover, the $(z,\xi)$-update gives the $z$-update \eqref{eq:QP:admm:z:new} in \texttt{ADMMSlack}:
\begin{lemma}[The $z$-update of \texttt{ADMMSlack} is the $(z,\xi)$- update in \eqref{eq:QP:zxi}]\label{lemma}
The solution to the QP \eqref{eq:QP:zxi} satisfies
$z =\widetilde{\mathrm{Proj}}_{[\ell,u]}(Ax+\mu/\rho)$
where \eqref{eq:proj:tilde} defines $\widetilde{\mathrm{Proj}}_{[\ell,u]}(\cdot)$.
\end{lemma}
Thus, \texttt{ADMMSlack} is a valid ADMM scheme for solving \QPslack, even though it does not  require the introduction of the slack variables $\xi$.

\begin{proof}
First, we consider separately the equality constraints $\ell_i=u_i=h_i$ for $i\in\mathcal{E}$ and the inequality constraints $\ell_i<u_i$ for $i\in\mathcal{I}$, with $\mathcal{E}\cup\mathcal{I}=\{1,\dots,p\}$. We then introduce the KKT  multipliers $\nu^h\in\R^{|\mathcal{E}|}$ for the equality constraints and $\nu^\ell,\nu^u\in\R^{|\mathcal{I}|}$ for the inequality constraints, and the Lagrangian of  the QP \eqref{eq:QP:zxi}
\begin{align}
\mathcal{L}(z,\xi,\nu)&=\tfrac{\rho}{2}\|z\|_2^2+\tfrac{\alpha}{2}\|\xi\|_2^2+(-\rho Ax-\mu)^\top z 
\\[1mm]
&\quad+ \sum_{i\in\mathcal{E}}\nu^h_i(\ell_i-z_i-\xi_i)
\nonumber
\\
&\quad+ \sum_{i\in\mathcal{I}}\nu^\ell_i(\ell_i-z_i-\xi_i)+\nu^u_i(z_i+\xi_i-u_i),
\nonumber
\end{align}
with a slight abuse of notation for the indices of the multipliers $\nu^h,\nu^\ell,\nu^u$. 

The KKT conditions associated to the QP \eqref{eq:QP:zxi} are:
\begin{subequations}
\label{eq:kkt}
\begin{align}
&\begin{cases}
\rho z_i-\rho A_ix-\mu_i -\nu^h_i = 0,
\\ 
\alpha \xi_i - \nu^h_i=0,
\\
\ell_i=z_i+\xi_i= u_i,
\end{cases}
&&i\in\mathcal{E},
\\
&\begin{cases}
\rho z_i-\rho A_ix-\mu_i -\nu^\ell_i+\nu^u_i = 0,
\\
\alpha \xi_i -\nu^\ell_i + \nu^u_i=0,
\\
\ell_i\leq z_i+\xi_i\leq u_i,
\\
\nu^\ell_i,\nu^u_i \geq 0,
\\
\nu^\ell_i(\ell_i-z_i-\xi_i)=0,
\\
\nu^u_i(z_i+\xi_i-u_i)=0,
\end{cases}
&&i\in\mathcal{I}.
\end{align}
\end{subequations}
Next, we study each index $i=1,\dots,p$ separately.

\textbf{1) Equality constraints, case  $\ell_i= u_i$, $i\in\mathcal{E}$.}\\
From the KKT conditions \eqref{eq:kkt}
\begin{align*}
z_i&=A_ix +\tfrac{1}{\rho}\mu_i+\tfrac{1}{\rho}\nu_i^h
\ \  \text{and} \ \   
\xi_i=\tfrac{1}{\alpha}\nu^h_i.
\end{align*}
We have the following two cases:

\textbf{1.1)} First, if $A_ix +\tfrac{1}{\rho}\mu_i\neq\ell_i$, then $\ell_i=z_i+\xi_i=A_ix +\tfrac{1}{\rho}\mu_i+\tfrac{\rho+\alpha}{\rho\alpha}\nu^h_i$, so $\nu^h_i=\tfrac{\rho\alpha}{\rho+\alpha}(\ell_i-A_ix-\tfrac{1}{\rho}\mu_i)$ and 
\begin{align*}
z_i&=A_ix+\tfrac{1}{\rho}\mu_i +\tfrac{\alpha}{\rho+\alpha}(\ell_i-
A_ix-\tfrac{1}{\rho}\mu_i)
\\
&=\frac{\rho A_ix+\mu_i+\alpha\ell_i}{\rho+\alpha}
\\
\xi_i &=\tfrac{\rho}{\rho+\alpha}(\ell_i-A_ix-\tfrac{1}{\rho}\mu_i).
\end{align*}

\textbf{1.2)} Second, if $A_ix +\tfrac{1}{\rho}\mu_i=\ell_i$, then $\ell_i=z_i+\xi_i=A_ix +\tfrac{1}{\rho}\mu_i+\tfrac{\rho+\alpha}{\rho\alpha}\nu^h_i$, so $\nu_i^h=0$ and 
\begin{align*}
z_i&=A_ix +\tfrac{1}{\rho}\mu_i 
\ \  \text{and} \ \   
\xi_i=0.
\end{align*}

\textbf{2) Inequality constraints, case  $\ell_i< u_i$, $i\in\mathcal{I}$.}\\
From the KKT conditions \eqref{eq:kkt},
\begin{align}
\label{eq:kkt:z}
&z_i=A_ix +\tfrac{1}{\rho}\mu_i+\tfrac{1}{\rho}(\nu^\ell_i-\nu^u_i),
\\
\label{eq:kkt:xi}
&\xi_i=\tfrac{1}{\alpha}(\nu^\ell_i - \nu^u_i),
\end{align}
we obtain 
\begin{align}
\label{eq:kkt:z+xi}
z_i+\xi_i=
A_ix+\tfrac{1}{\rho}\mu_i +\tfrac{\rho+\alpha}{\rho\alpha}(\nu^\ell_i - \nu^u_i),
\\
\label{eq:nul-nuu}
\nu^\ell_i - \nu^u_i
=
\tfrac{\rho\alpha}{\rho+\alpha}(z_i+\xi_i-
A_ix-\tfrac{1}{\rho}\mu_i).
\end{align}
Also, by \eqref{eq:kkt} and $\ell_i\neq u_i$, we have that
\begin{align}\label{eq:nuu>0}
\nu^u_i>0 \implies \nu^\ell_i=0,
\\
\label{eq:nuell>0}
\nu^\ell_i>0 \implies \nu^u_i=0.
\end{align}
Then, we have the following three cases:

\textbf{2.1)} First, if $A_ix+\tfrac{1}{\rho}\mu_i<\ell_i$, then %
\eqref{eq:kkt:z+xi} implies $z_i+\xi_i<\ell_i+\tfrac{\rho+\alpha}{\rho\alpha}(\nu^\ell_i - \nu^u_i)$. Also, $\ell_i\leq z_i+\xi_i$, so $$\ell_i<\ell_i+\tfrac{\rho+\alpha}{\rho\alpha}(\nu^\ell - \nu^u).
$$
\begin{itemize}
\item If $\nu^u_i>0$, then $\nu^\ell_i=0$   by \eqref{eq:nuu>0}. But in this case, we arrive at a contradiction with $\ell_i<\ell_i+\tfrac{\rho+\alpha}{\rho\alpha}(- \nu^u)$.
\item If $\nu^u_i=0$ and $\nu^\ell_i=0$, then we also  arrive at a contradiction with $\ell_i<\ell_i$.
\item Thus, we conclude that $\nu^\ell_i>0$  and $\nu^u_i=0$.
\end{itemize}
Since  $\nu_i^\ell>0$, we have $z_i+\xi_i=\ell_i$ by \eqref{eq:kkt}. 
Then, \eqref{eq:nul-nuu} implies that $\nu_i^\ell =\frac{\rho\alpha}{\rho+\alpha}(\ell_i-
A_ix-\tfrac{1}{\rho}\mu_i)$, so 
\begin{align*}
z_i&\mathop{=}^\eqref{eq:kkt:z}
A_ix+\tfrac{1}{\rho}\mu_i +\tfrac{\alpha}{\rho+\alpha}(\ell_i-
A_ix-\tfrac{1}{\rho}\mu_i)
\\
&=\frac{\rho A_ix+\mu_i+\alpha\ell_i}{\rho+\alpha},
\\
\xi_i &\mathop{=}^\eqref{eq:kkt:xi}\tfrac{\rho}{\rho+\alpha}(\ell_i-A_ix-\tfrac{1}{\rho}\mu_i).
\end{align*}

\textbf{2.2)} Second, if $A_ix+\tfrac{1}{\rho}\mu_i>u_i$, similar arguments as described above imply
\begin{align*}
z_i
&=\frac{\rho A_ix+\mu_i+\alpha u_i}{\rho+\alpha},
\\
\xi_i &=\tfrac{\rho}{\rho+\alpha}(u_i-A_ix-\tfrac{1}{\rho}\mu_i).
\end{align*}

\textbf{2.3)} 
Third, if $\ell_i\leq A_ix+\tfrac{1}{\rho}\mu_i\leq u_i$, then we look at three cases:
\begin{itemize}
\item $\nu^\ell_i>0$ and $\nu^u_i>0$ leads to a contradiction with \eqref{eq:kkt} since $\ell_i\neq u_i$. 
\item If $\nu^\ell_i=0$ and $\nu^u_i\geq 0$, then:
\begin{itemize}
\item If $\nu^u_i>0$, then $z_i+\xi_i=u_i$ by \eqref{eq:kkt}, so 
\begin{align*}
u_i&\mathop{=}^\eqref{eq:kkt:z+xi} A_ix+\tfrac{1}{\rho}\mu_i -\tfrac{\rho+\alpha}{\rho\alpha}\nu^u_i\leq u_i-\tfrac{\rho+\alpha}{\rho\alpha}\nu^u_i<u_i,
\end{align*}
which is a contradiction.
\item Thus, $\nu^u_i=0$, and we obtain
\end{itemize} 
\begin{align*}
z_i&\,\mathop{=}^\eqref{eq:kkt:z}\,A_ix+\tfrac{1}{\rho}\mu_i, 
\ \, 
\xi_i\,\mathop{=}^\eqref{eq:kkt:xi}\,0.
\end{align*}
\item If $\nu^\ell_i=0$ and $\nu^u_i\geq 0$, then similar arguments show that $\nu^u_i=0$, so the equations above hold. 
\end{itemize}
To conclude,  the solution to the QP \eqref{eq:QP:zxi} is
\begin{align*}
&\begin{cases}
\begin{cases}
z_i=\frac{\rho A_ix+\mu_i+\alpha\ell_i}{\rho+\alpha}
\\ 
\xi_i =\tfrac{\rho}{\rho+\alpha}(\ell_i-A_ix-\tfrac{1}{\rho}\mu_i),
\end{cases}
&\hspace{-5mm}\text{if }A_ix +\tfrac{1}{\rho}\mu_i\neq\ell_i,
\\ 
\begin{cases}
z_i=A_ix +\tfrac{1}{\rho}\mu_i,
\\
\xi_i=0
\end{cases}
&\hspace{-5mm}\text{if }A_ix +\tfrac{1}{\rho}\mu_i=\ell_i,
\end{cases}
\end{align*}
for $i\in\mathcal{E}$, and
\begin{align*}
&\begin{cases}
\begin{cases}
z_i=\frac{\rho A_ix+\mu_i+\alpha\ell_i}{\rho+\alpha}
\\ 
\xi_i =\tfrac{\rho}{\rho+\alpha}(\ell_i-A_ix-\tfrac{1}{\rho}\mu_i),
\end{cases}
&\hspace{-4mm}\text{if }A_ix +\tfrac{1}{\rho}\mu_i<\ell_i,
\\ 
\begin{cases}
z_i=\frac{\rho A_ix+\mu_i+\alpha u_i}{\rho+\alpha}
\\ 
\xi_i =\tfrac{\rho}{\rho+\alpha}(u_i-A_ix-\tfrac{1}{\rho}\mu_i),
\end{cases}
&\hspace{-4mm}\text{if }A_ix +\tfrac{1}{\rho}\mu_i> u_i,
\\ 
\begin{cases}
z_i=A_ix +\tfrac{1}{\rho}\mu_i,
\\
\xi_i=0
\end{cases}
&\hspace{-11mm}\text{if }\ell_i\leq A_ix +\tfrac{1}{\rho}\mu_i\leq u_i,
\end{cases}
\end{align*}
for $i\in\mathcal{I}$. Thus, we obtain   that $z=\widetilde{\text{Proj}}_{[\ell,u]}(Ax+\mu/\rho)$, which concludes the proof of Lemma \ref{lemma}.
\end{proof}

\begin{remark}[Interpretation as an averaged projection]\label{remark:averagedprojection}
The $z$-update of \texttt{ADMMSlack} is equivalently  written as\hspace{-3mm}
$$
z = \left(1-\gamma\right)\tilde{z} + \gamma\,\text{Proj}_{[\ell,u]}(\tilde{z})=\eqref{eq:QP:admm:z:new},
$$
where $\gamma:=\frac{\alpha}{\alpha+\rho}$ and $\tilde{z}:=Ax+\mu/\rho$. As $\alpha$ increases, the $z$-update becomes a hard projection. As $\alpha$ decreases, the $z$-update becomes an averaged projection. Averaged projections of this form are standard in the proximal operators literature \cite[Section 2.3]{Parikh2014} \cite[Section 3]{Eckstein1992}.
\end{remark}

\section{Numerical results}\label{sec:results}
We compare the proposed \texttt{ADMMSlack} scheme  with the standard \texttt{ADMM} scheme  for problems with slack variables, resulting in the algorithms summarized in Table \ref{tab:admm}.
The two ADMM schemes are evaluated on random
MPC problems with slack variables:
\begin{subequations}
\label{eq:trajopt}
\begin{align}
\hspace{-2mm}
\min_{\substack{x,u,\xi}}\ \   
&\sum_{k=0}^N \tfrac{1}{2}\|x_k\|^2_Q + \sum_{k=0}^{N-1}\tfrac{1}{2}\|u_k\|^2_R + \tfrac{\alpha}{2}\|\xi\|^2
\nonumber
\\
\label{eq:trajopt:dynamics}
\hspace{-2mm}
\textrm{s.t.} \quad
&x_{k+1}=Ax_k+Bu_k+\xi_k^A,
\ \   k=0,\dots,N,
\\
&x_0 = x_s + \xi^s, 
\\
\label{eq:trajopt:constraints:state}
&\underline{x}\leq x_k+\xi^x_k\leq \overline{x},
\qquad\ \ \ k=0,\dots,N,
\\ 
\label{eq:trajopt:constraints:control}
&\underline{u}\leq u_k+\xi^u_k\leq \overline{u},
\qquad\ \ \ k=0,\dots,N{-}1,
\end{align}
\end{subequations}
with $\|x\|^2_Q:=x^\top Q x$, %
$N=20$, $\alpha=10$, and the slack variables $\xi:=((\xi^A_k)_{k=0}^N,\xi^s,(\xi^x_k)_{k=0}^N,(\xi^u_k)_{k=0}^{N-1})$. 

We randomly sample the MPC problem parameters $(Q,R,A,B,\underline{x},\overline{x},\underline{u},\overline{u})$ %
as in \cite{Stellato2020}. 
The initial states are sampled according to two scenarios:
\begin{itemize}
    \item \textbf{Feasible QPs}: We sample the initial states $x_s$ strictly inside the box constraints $[\underline{x},\overline{x}]$ as in \cite{Stellato2020}.
    \item \textbf{Infeasible QPs}: We sample the initial states  $x_s$ strictly outside the box constraints $[\underline{x},\overline{x}]$ as $x_s\sim (1+w)\overline{x}$ with $w$ uniformly-distributed in $[0,1]$. These MPC problems would be infeasible without the slack variables $\xi$, and they correspond to MPC applications where the system starts from initial conditions that slightly violate constraints, e.g., due to model mismatch in real-world deployment.   
\end{itemize}
\textbf{Code}: The ADMM solvers are implemented in  \texttt{C++} and \texttt{Python} for both a performant and easily readable implementation. %
Computation times are evaluated using the \texttt{C++} implementation on a laptop with a 12th Gen Intel Core i7-1280Px20 CPU with 32 GB of RAM. 
The code to reproduce results is available at:
\\ 
\texttt{\url{https://github.com/A2R-Lab/ADMMSlack}}

\subsection{\texttt{ADMMSlack} has a similar convergence behavior}

First, we show that \texttt{ADMM}  and  \texttt{ADMMSlack}  have similar convergence behavior, although \texttt{ADMMSlack} solves smaller linear systems in the $x-$updates (see Table \ref{tab:admm}).  

\begin{figure}[!t]
    \centering
    \textbf{Feasible QPs}
    \\[1mm]
    \includegraphics[width=\linewidth]{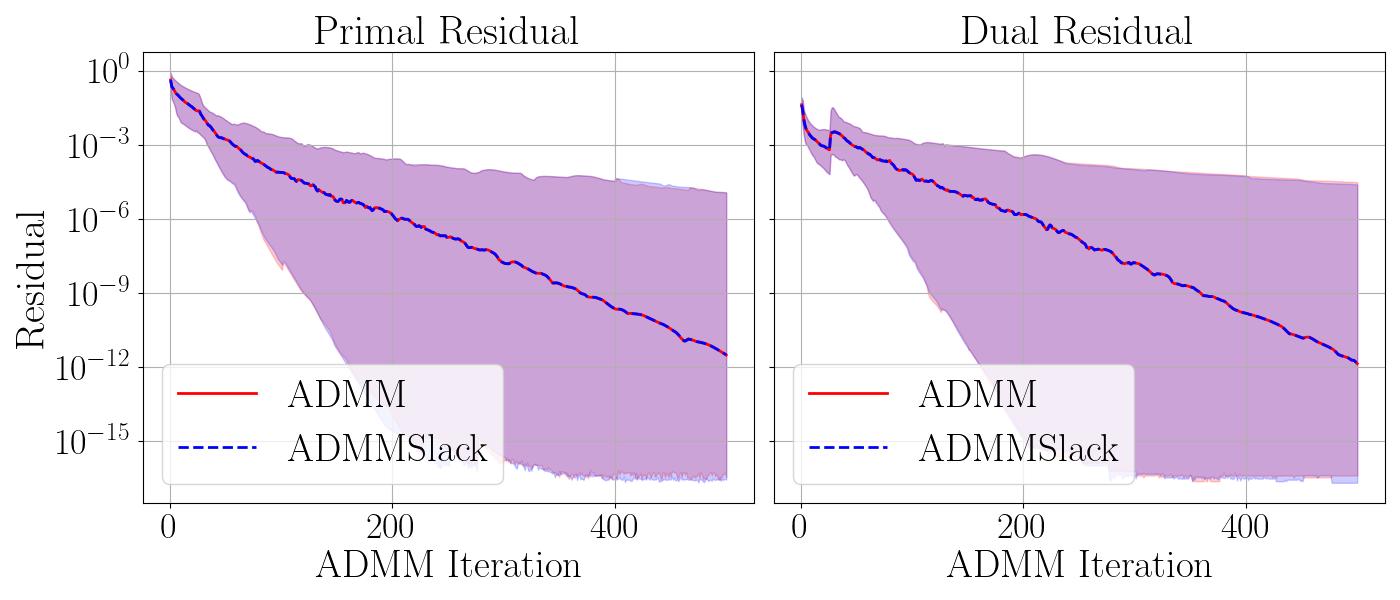}
    \\
    \textbf{Infeasible QPs}
    \\[1mm]
    \includegraphics[width=\linewidth]{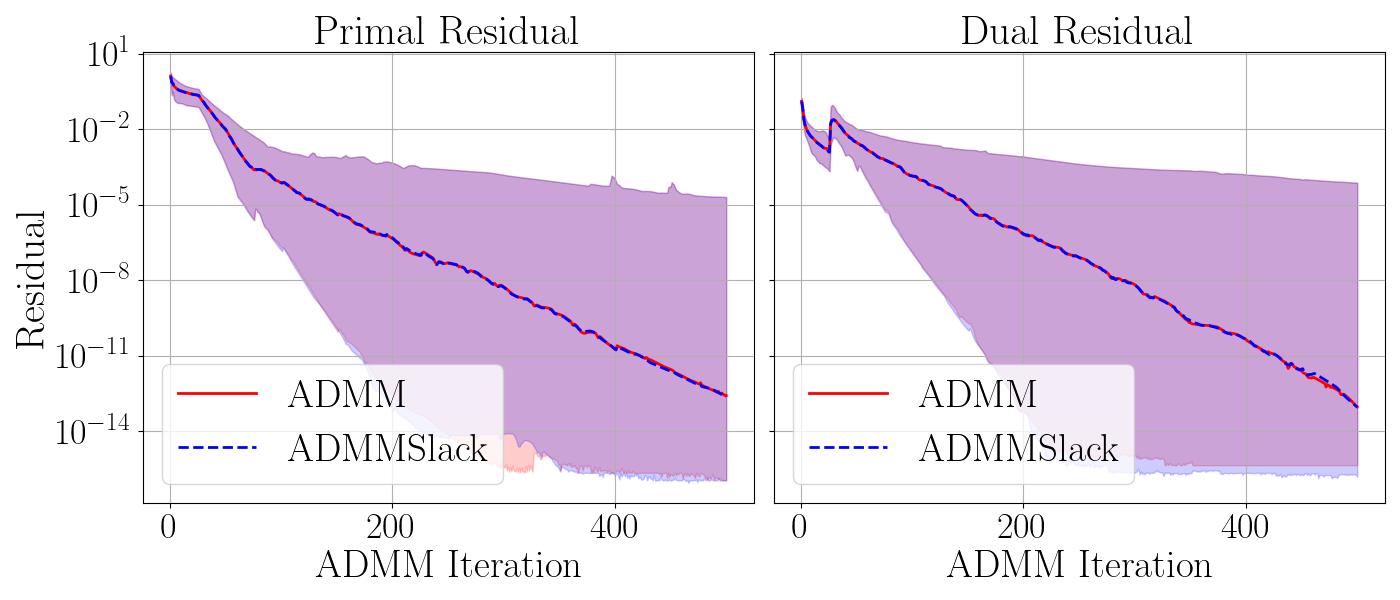}
    \caption{Residuals across ADMM iterates for the feasible (top) and infeasible QPs (bottom): median over $100$ MPC problem instances, with $\pm 2$ standard deviations.}
    \label{fig:results:residuals}
\end{figure}
\begin{figure}[!t]
    \centering
    \textbf{Feasible QPs}
    \\[1mm]
    \includegraphics[width=\linewidth]{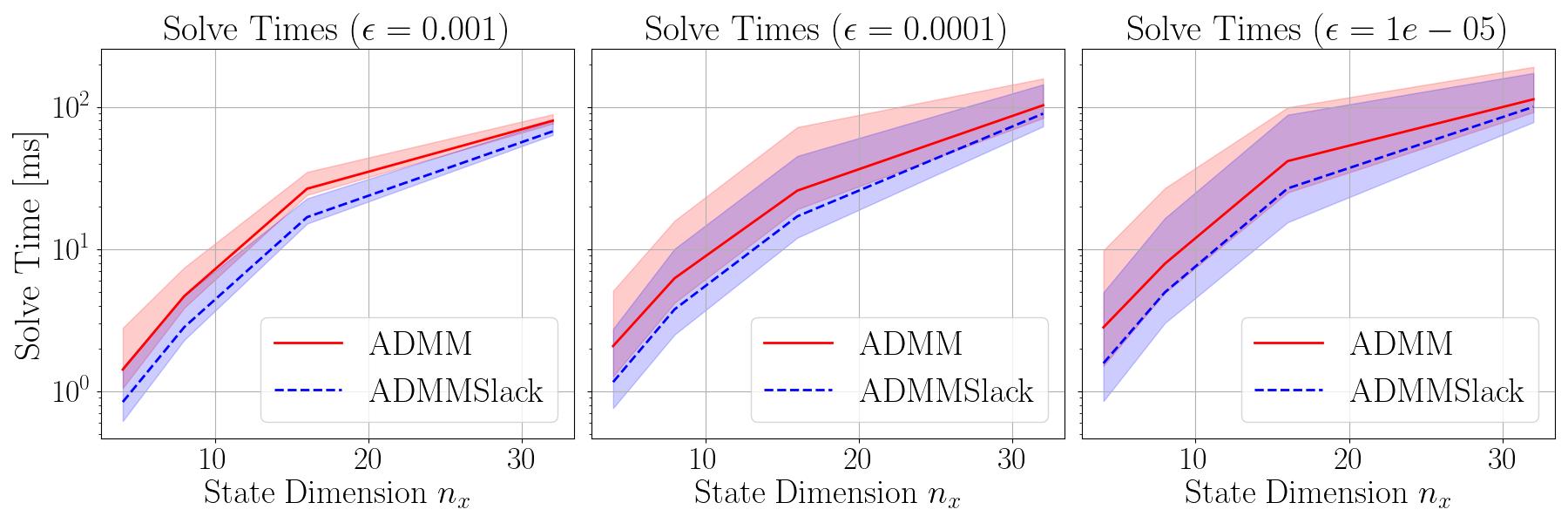}
    \\
    \hspace{0.01\linewidth}\includegraphics[width=0.99\linewidth]{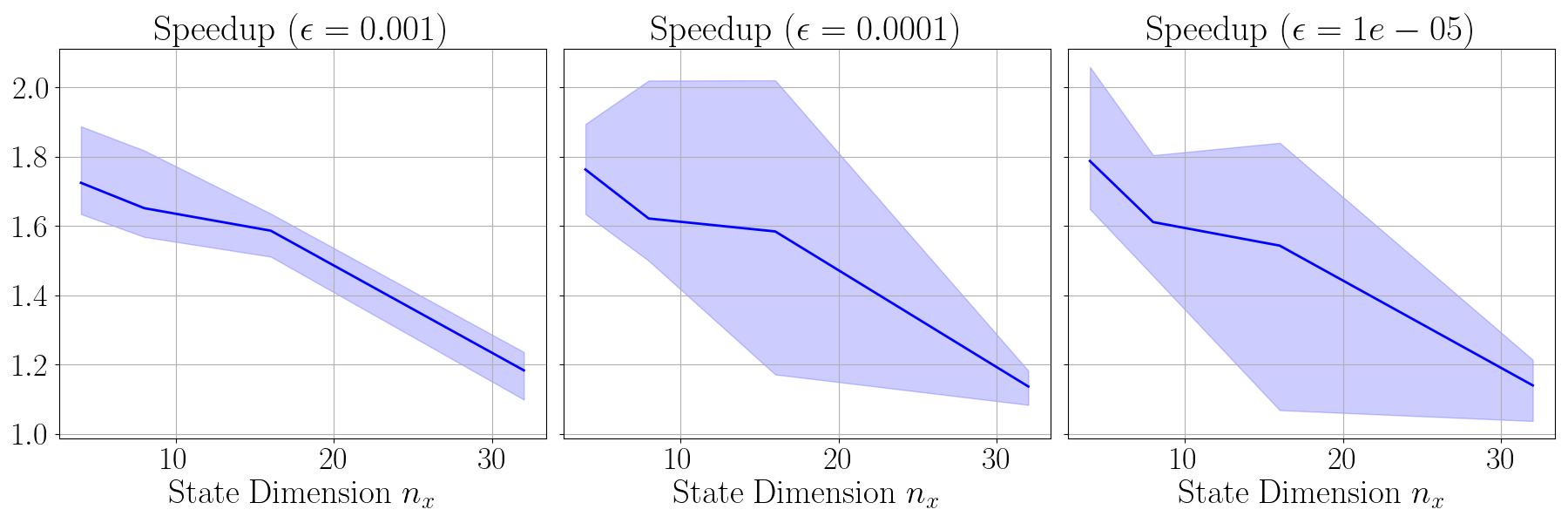}
    \\
    \textbf{Infeasible QPs}
    \\[1mm]
    \includegraphics[width=\linewidth]{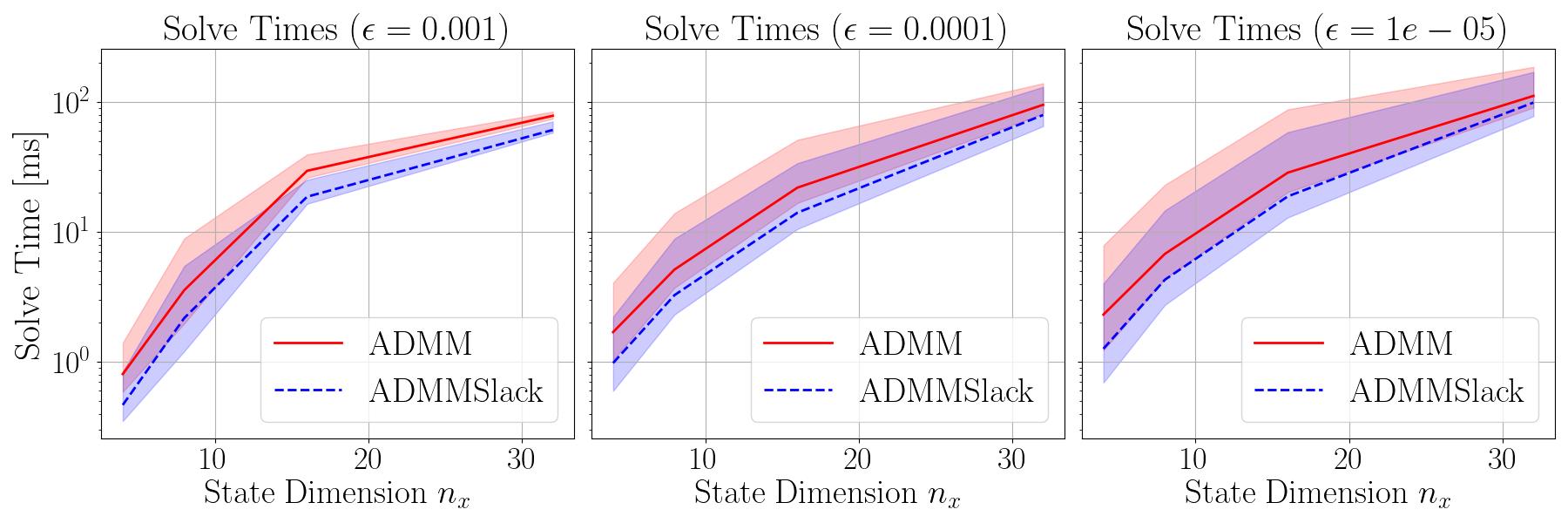}
    \\
    \hspace{0.01\linewidth}\includegraphics[width=0.99\linewidth]{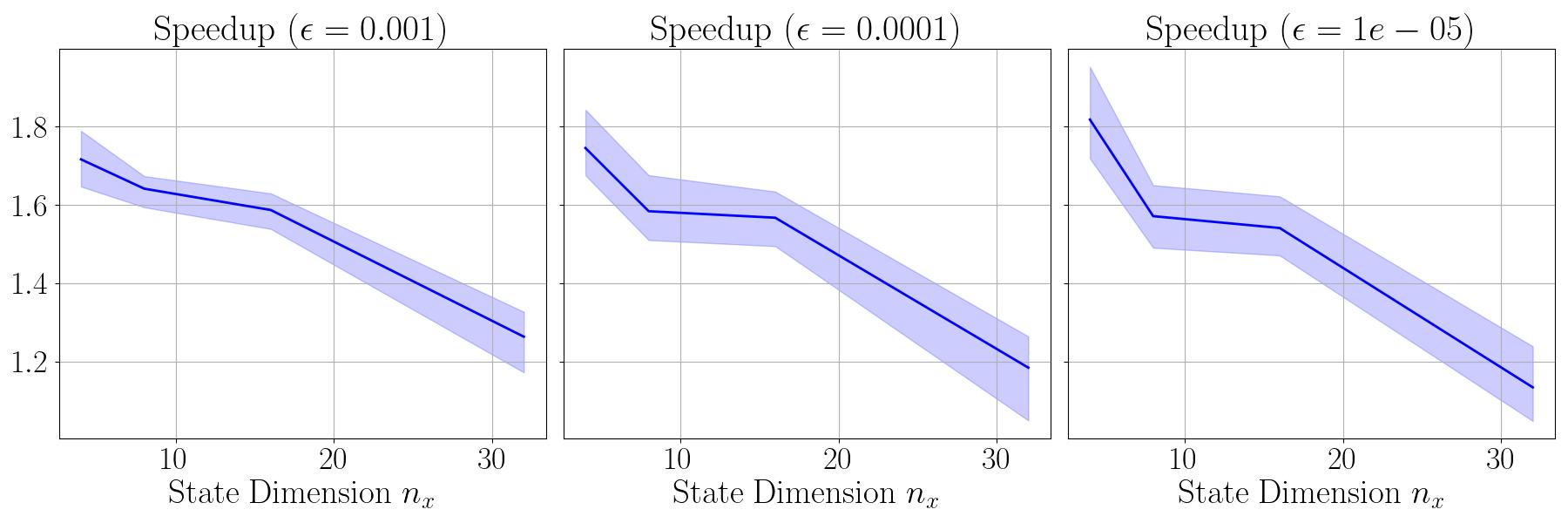}
    \caption{Solve times and speedups (ratio $\big(\text{SolveTime}(\texttt{ADMM})$ / $\text{SolveTime}(\texttt{ADMMSlack})\big)$) for the feasible (top) and infeasible QPs (bottom) for different exit tolerances and state-control dimensions: median over $100$ MPC problem instances with $\pm 2$ standard deviations.
    \vspace{-2mm}
    }
    \label{fig:results:speed}
\end{figure}

We set $\text{dim}(x_k)=2\dim(u_k)=4$ and solve $100$ random MPC problem instances with the two ADMM schemes. We report the primal and dual residuals across ADMM iterations in Figure~\ref{fig:results:residuals}. \texttt{ADMMSlack} exhibits similar convergence behavior as \texttt{ADMM}, as both the median and $\pm 2$ standard deviations of the residuals remain mostly identical across both feasible and infeasible QPs. 
Overall, these results show that the two schemes take a similar  number of iterations to reach a desired exit tolerance for these problems.

\subsection{\texttt{ADMMSlack} is faster}
We vary the exit tolerance $\epsilon$ and the state and control dimensions $\text{dim}(x_k)=2\dim(u_k)$, and for each variation, again solve $100$ random MPC problem instances with the two ADMM schemes. In Figure~\ref{fig:results:speed} we report the median and $\pm 2$ standard deviations of the resulting solve times.
Results show that the proposed \texttt{ADMMSlack} scheme converges significantly faster than the standard \texttt{ADMM} scheme over all problem instances, for both feasible and infeasible QPs and for different state-control dimensions and exit tolerances $\epsilon$. This speedup is due to solving smaller linear systems in $x$-update, as these linear systems do not involve the slack variables $\xi$.

The performance of  \texttt{ADMM} and  \texttt{ADMMSlack}  depends on problem conditioning (affected by $\alpha$) and the step size parameters $\rho$. For fair comparisons, we used identical parameters for both methods. We also varied $\alpha\in\{0.1,1,10^2,10^3\}$ and obtained similar results  as reported above. %

\section{Conclusion}
Handling slack variables with an ADMM scheme does not necessarily require increasing the size of the problem. This finding could be integrated into other solvers to improve their performance for solving quadratic programs with slack variables.

While we focus on quadratic penalization of the slack variables $\|\xi\|^2$,  it would be interesting to investigate other penalizations such as the $1$-norm  $\|\xi\|_1$.

\bibliographystyle{IEEEtran}
\bibliography{main}
\end{document}

%% file: preamble.tex
\usepackage{wrapfig}

\usepackage{graphicx}
\usepackage{amsmath}
\usepackage{amssymb}
\usepackage{bm}
\usepackage{accents}

\newcommand{\uwidehat}[1]{%
  \mathpalette\douwidehat{#1}%
}
\makeatletter
\newcommand{\douwidehat}[2]{%
  \sbox0{$\m@th#1\widehat{\hphantom{#2}}$}%
  \sbox2{$\m@th#1x$}
  \sbox4{$\m@th#1#2$}
  \dimen0=\ht0
  \advance\dimen0 -.8\ht2
  \dimen2=\dp4
  \rlap{%
    \raisebox{\dimexpr\dimen0-\dimen2}{%
      \scalebox{1}[-1]{\box0}%
    }%
  }%
  {#2}%
}
\makeatother

\def\namedlabel#1#2{\begingroup
    #2%
    \def\@currentlabel{#2}%
    \phantomsection\label{#1}\endgroup
}
\makeatletter
\newcommand{\leqnomode}{\tagsleft@true\let\veqno\@@leqno}
\makeatother

\makeatletter
\newcommand{\customlabel}[2]{%
\protected@write\@auxout{}{\string \newlabel{#1}{{#2}{\thepage}{#2}{#1}{}}}%
\hypertarget{#1}{#2}
}
\makeatother

\newcommand\PMPODE[1]{\hyperref[PMPODE]{$\textbf{ODE}_{#1}$}\xspace}
\newcommand\PMPODElambda[1]{\hyperref[PMPODElambda]{$\textbf{ODE}_{#1}^\lambda$}\xspace}
\newcommand\PMPODElambdaunder[1]{\hyperref[PMPODElambda_under]{$\uwidehat{\textbf{ODE}_{#1}^\lambda}$}\xspace}
\newcommand\PMPODElambdaover[1]{\hyperref[PMPODElambda_over]{$\widehat{\textbf{ODE}_{#1}^\lambda}$}\xspace}
\newcommand\PMPODEepsilonfullrank[1]{\hyperref[PMPODEepsilon_full_rank]{$\textbf{ODE}_{#1}^\epsilon$}\xspace}

\newcommand\PMPODEg[1]{\hyperref[PMPODEg]{$\textbf{ODE}^g_{#1}$}\xspace}
\newcommand\OCPg[1]{\hyperref[OCPg]{$\textbf{OCP}^g_{#1}$}\xspace}

\newcommand\OCP[1]{\hyperref[OCP]{$\textbf{OCP}_{#1}$}\xspace}
\newcommand\BVP[1]{\hyperref[BVP]{$\textbf{BVP}_{#1}$}\xspace}
\newcommand\OCPMPC[1]{\hyperref[OCPMPC]{$\textbf{OCP}(#1)$}\xspace}

\usepackage{tabularx}
\newcolumntype{C}{>{\centering\arraybackslash}p{19mm}}
\newcolumntype{G}{>{\centering\arraybackslash}p{4mm}}
\newcolumntype{S}{>{\centering\arraybackslash\scriptsize}p{4mm}}

\usepackage[font=small,labelfont=bf]{caption}
\usepackage{color}
\usepackage[usenames,dvipsnames]{xcolor}

\usepackage{cleveref} %

\usepackage{xspace}

\let\labelindent\relax
\usepackage{enumitem} %

\usepackage{cite} %

\usepackage{amssymb}%
\usepackage{algorithm}
\usepackage{algcompatible}
\makeatletter
\renewcommand*{\ALG@name}{Alg.}
\makeatother
\usepackage[noend]{algpseudocode}

\usepackage{xpatch} %

\usepackage{tabularx, booktabs} %
\newcolumntype{F}{>{\centering\arraybackslash}p{1.9cm}}

\newcolumntype{G}{>{\centering\arraybackslash}p{1.2cm}}
\newcolumntype{H}{>{\centering\arraybackslash}p{0.4cm}}
\newcolumntype{I}{>{\centering\arraybackslash}p{0.92cm}}
\newcolumntype{Z}{>{\centering\arraybackslash}p{-0.92cm}}

\newtheorem{assumption}{Assumption}
\newtheorem{lemma}{Lemma}

\newtheorem{remark}{Remark}

\newcommand\mydots{\hbox to 1em{.\hss.\hss.}}

\newcommand{\R}{\mathbb{R}}

\newtheorem{preremark3}{Theorem}[section]

\usepackage{mdframed}
\usepackage{lipsum}
\newmdtheoremenv{theo}{Theorem}

\usepackage{multirow}

\usepackage{accents}

\newcommand\bigO[1]{{\color{black!50!white}O(#1)}}

%% file: table.tex
\begin{table*}[t!]
\centering
\renewcommand{\arraystretch}{1.6}
\begin{tabular}{c|p{4.6cm}||p{5.1cm}|p{4.4cm}}
\toprule
\hspace{-1mm}\textbf{Update}\hspace{-1mm}
 & 
\parbox[c]{4.6cm}{
\rule{0pt}{4mm}
\centering\textbf{ADMM for \QP}
\\
Solves the QP without slack variables
}
& 
\parbox[c]{4.6cm}{
\rule{0pt}{4mm}
\centering\textbf{ADMM for \QPslack} 
\\
Introduces additional slack variables
}
&
\parbox[c]{4.2cm}{
\rule{0pt}{4mm}
\centering \texttt{ADMMSlack} for \QPslack
\\
Does not introduce slack variables}
\\[3mm]
\midrule
\textbf{$\mathbf{x}$} 
& \eqref{eq:QP:admm:x} $\, x = (Q + \rho A^\top A)^{-1}(\cdot )$ \hfill $\bigO{n^3}$
& \eqref{eq:QP:admm:x:slack} $\, \bar{x} = (\bar{Q} + \rho \bar{A}^\top \bar{A})^{-1}(\cdot)$ \hfill $\bigO{(n+p)^3}$
& \eqref{eq:QP:admm:x} $\, x = (Q + \rho A^\top A)^{-1}(\cdot)$ \hfill $\bigO{n^3}$
\\
\textbf{$\mathbf{z}$} 
& \eqref{eq:QP:admm:z} $\, z = \text{Proj}_{[\ell, u]}(Ax + \mu/\rho)$ \hfill $\bigO{np}$
&  \eqref{eq:QP:admm:z:slack} $\, z = \text{Proj}_{[\ell, u]}(\bar{A}\bar{x} + \mu/\rho)$  \hfill $\bigO{np}$
&  \eqref{eq:proj:tilde} $\, z = \widetilde{\text{Proj}}_{[\ell, u]}(Ax + \mu/\rho)$ \hfill  $\bigO{np}$
\\
\textbf{dual} 
&  \eqref{eq:QP:admm:mu} $\, \mu = \mu + \rho(Ax - z)$  \hfill $\bigO{np}$
&  \eqref{eq:QP:admm:mu:slack} $\, \mu = \mu + \rho(\bar{A}\bar{x} - z)$  \hfill $\bigO{np}$
&  \eqref{eq:QP:admm:mu} $\, \mu = \mu + \rho(Ax - z)$  \hfill $\bigO{np}$\\
\bottomrule
\end{tabular}
\caption{ADMM algorithms for solving \QP and \QPslack. Compared to a standard ADMM scheme applied to \QPslack, the proposed ADMM scheme (\texttt{ADMMSlack}) does not increase the problem size. The only difference between \texttt{ADMMSlack} and an ADMM scheme for \QP is a slight modification of the projection operation in the $z$-update step. %
\vspace{-3mm}
}
\label{tab:admm}
\end{table*}

%% file: main.bib
@String { jrn_IEEE_TAC              = {{IEEE Transactions on Automatic Control}} }

@String { proc_EUCA_ECC             = {{European Control Conference}} }

@String { proc_IEEE_ICRA            = {{Proc.\ IEEE Conf.\ on Robotics and Automation}} }

@incollection{Houska2018,
  year = {2018},
  publisher = {Springer International Publishing},
  pages = {413--443},
  author = {Houska, B. and Villanueva, M.~E.},
  title = {Robust Optimization for {MPC}},
  booktitle = {Handbook of Model Predictive Control}
}

@article{Stellato2020,
  year = {2020},
  volume = {12},
  number = {4},
  pages = {637--672},
  author = {Stellato, B. and Banjac, G. and Goulart, P. and Bemporad, A. and Boyd, S.},
  title = {{OSQP}: an operator splitting solver for quadratic programs},
  journal = {Mathematical Programming Computation}
}

@inproceedings{adabag2024mpcgpu,
  title={{MPCGPU}: Real-Time Nonlinear Model Predictive Control through Preconditioned Conjugate Gradient on the {GPU}}, 
  author={Adabag, E. and Atal, M. and Gerard, W. and Plancher, B.},
  booktitle=proc_IEEE_ICRA,
  year = {2024}
}

@inproceedings{LewGreiffICRA2025,
  title = {Risk-Averse Model Predictive Control for Racing in Adverse Conditions},
  booktitle = proc_IEEE_ICRA,
  author       = {Lew, T. and Greiff, M. and Djeumou, F. and Suminaka, M. and Thompson, M. and Subosits, J.},
  year = {2025},
}

@article{Boyd2011,
  title = {Distributed Optimization and Statistical Learning via the Alternating Direction Method of Multipliers},
  volume = {3},
  number = {1},
  journal = {Foundations and Trends{\textregistered} in Machine Learning},
  publisher = {Foundations and Trends in Machine Learning},
  author = {Boyd, S. and Parikh, N. and Chu, E. and Peleato, B. and Eckstein, J.},
  year = {2011},
  pages = {1--122}
}

@article{venkat2008distributed,
  title={Distributed {MPC} strategies with application to power system automatic generation control},
  author={Venkat, Aswin N and Hiskens, Ian A and Rawlings, James B and Wright, Stephen J},
  journal={IEEE transactions on control systems technology},
  volume={16},
  number={6},
  pages={1192--1206},
  year={2008},
  publisher={IEEE}
}

@article{o2016conic,
  title={Conic optimization via operator splitting and homogeneous self-dual embedding},
  author={O’donoghue, Brendan and Chu, Eric and Parikh, Neal and Boyd, Stephen},
  journal={Journal of Optimization Theory and Applications},
  volume={169},
  number={3},
  pages={1042--1068},
  year={2016},
  publisher={Springer}
}

@book{del2010automotive,
  title={Automotive model predictive control},
  author={Del Re, Luigi and Allgower, Frank and Glielmo, Luigi and Guardiola, Carlos and Kolmanovsky, Ilya},
  year={2010},
  publisher={Springer}
}

@article{tranzatto2022cerberus,
  title={Cerberus in the {DARPA} {S}ubterranean {C}hallenge},
  author={Tranzatto, Marco and Miki, Takahiro and Dharmadhikari, Mihir and Bernreiter, Lukas and Kulkarni, Mihir and Mascarich, Frank and Andersson, Olov and Khattak, Shehryar and Hutter, Marco and Siegwart, Roland and others},
  journal={Science Robotics},
  volume={7},
  number={66},
  pages={eabp9742},
  year={2022},
  publisher={American Association for the Advancement of Science}
}

@article{wensing2023optimization,
  title={Optimization-based control for dynamic legged robots},
  author={Wensing, Patrick M and Posa, Michael and Hu, Yue and Escande, Adrien and Mansard, Nicolas and Del Prete, Andrea},
  journal={IEEE Transactions on Robotics},
  year={2023},
  publisher={IEEE}
}

@article{vazquez2014model,
  title={Model predictive control: A review of its applications in power electronics},
  author={Vazquez, Sergio and Leon, Jose I and Franquelo, Leopoldo G and Rodriguez, Jose and Young, Hector A and Marquez, Abraham and Zanchetta, Pericle},
  journal={IEEE industrial electronics magazine},
  volume={8},
  number={1},
  pages={16--31},
  year={2014},
  publisher={IEEE}
}

@article{Feller2017,
  title = {Relaxed Logarithmic Barrier Function Based Model Predictive Control of Linear Systems},
  volume = {62},
  number = {3},
  journal = jrn_IEEE_TAC,
  author = {Feller,  Christian and Ebenbauer,  Christian},
  year = {2017},
  pages = {1223--1238}
}

@inproceedings{Oancea2023,
  title = {Relaxed Recentered Log-Barrier Function Based Nonlinear Model Predictive Control},
  booktitle = proc_EUCA_ECC,
  author = {Oancea, T.~A. and Jiang, Y. and Jones, C.~N.},
  year = {2023},
}

@article{Parikh2014,
  title = {Proximal Algorithms},
  volume = {1},
  number = {3},
  journal = {Foundations and Trends in Optimization},
  author = {Parikh, N. and Boyd, S.},
  year = {2014},
  pages = {127--239}
}

@article{Eckstein1992,
  title = {On the {Douglas-Rachford} splitting method and the proximal point algorithm for maximal monotone operators},
  volume = {55},
  number = {1-3},
  journal = {Mathematical Programming},
  author = {Eckstein,  J. and Bertsekas,  D.~P.},
  year = {1992},
  pages = {293--318}
}
